\documentclass{article}

 \usepackage[preprint]{neurips_2026}
 \newif\ifarxiv
 \arxivtrue

\usepackage[utf8]{inputenc} 
\usepackage[T1]{fontenc}    
\usepackage{hyperref}       
\usepackage{url}            
\usepackage{booktabs}       
\usepackage{amsfonts}       
\usepackage{nicefrac}       
\usepackage{microtype}      
\usepackage{xcolor}         

\usepackage{graphicx}
\usepackage{comment}
\usepackage{amsmath}
\usepackage{amssymb}
\usepackage{mathtools}
\usepackage{algorithm}      
\usepackage{algpseudocode}
\usepackage{enumitem}
\usepackage[bottom]{footmisc}
\usepackage{lipsum}

\newcommand{\inred}[1]{{\color{red}#1}}

\title{A Differentiable Interior-Point Method\\ in Single Precision}

\author{%
  Jon Arrizabalaga\\
  MIT\\
  Cambridge, MA 02139 \\
  \texttt{jonarri@mit.edu} \\
  \And
  Kevin Tracy\\
  Gridmatic\\
  Cupertino, CA 95014 \\
  \texttt{kevintracy48@gmail.com} \\
  \And
  Zachary Manchester\\
  MIT\\
  Cambridge, MA 02139 \\
  \texttt{zacm@mit.edu}
}

\begin{document}

\maketitle

\begin{abstract}

Primal-dual interior-point methods solve constrained convex optimization problems to tight tolerances with speed and robustness. Their solutions are also efficiently differentiable with respect to the problem data through the implicit function theorem. However, the standard treatment of primal-dual complementarity makes the underlying linear systems increasingly ill-conditioned near the solution. While this ill-conditioning is often benign in double precision, it can be catastrophic in single precision, preventing interior-point methods from fully exploiting the accelerated hardware that underpins modern machine learning. This paper introduces a differentiable interior-point method designed for low-precision arithmetic. By using an alternative complementarity representation, we ensure that the underlying linear systems remain spectrally bounded --- even near the solution --- a property that is essential for computing accurate gradients and avoiding arithmetic exceptions. As a result, our method enables interior-point solvers to reliably solve and differentiate optimization problems in single precision that were previously confined to double precision. We demonstrate the approach through an ablation study against the standard interior-point formulation and applications in bilevel and end-to-end learning settings where differentiating through constrained optimization is essential. The source code is available at \url{https://github.com/qpax-solver/qpax}.

\end{abstract}

\section{Introduction}~\label{sec:introduction}

Constrained convex programs are central to control~\cite{kuindersma2014efficiently, blackmore2016autonomous, varin2018constrained, 9905530, vega2025convexmaneuverplanningspacecraft}, simulation~\cite{anitescu2006optimization,tracy2023differentiable,pang2021convex}, operations research~\cite{beckmann1956studies, scarf1958minmax, ben1998robust,rubinstein2002markowitz,bertsimas2004price}, and machine learning~\cite{cortes1995support,tibshirani1996regression,candes2009exact,friedman2008sparse,boyd2011distributed}. They provide a principled way to encode domain-specific structure, often leading to optimization problems that can be solved efficiently and reliably~\cite{boyd2004convex}. With the growth of large-scale data and compute, increasing attention has shifted toward making these optimization pipelines differentiable, so that problem parameters can be tuned directly from data~\cite{donti2017task, wilder2019melding,salzmann24a}. While this does not replace expert modeling, it substantially expands the design space by allowing optimization-based components to be embedded within larger trainable systems~\cite{10049101,11301631,adabag2025differentiablemodelpredictivecontrol}.

In recent years, differentiable convex optimization has attracted considerable interest, particularly as a mechanism for integrating optimization layers into neural architectures~\cite{agrawal2019differentiating, amos2017optnet,agrawal2019differentiable,tracy2024differentiability}. This line of work has led to practical modeling tools, such as CVXPY Layers, that make it straightforward to incorporate parameterized convex programs into standard machine learning workflows~\cite{diamond2016cvxpy}.

Primal-dual interior-point methods (IPMs) have long served as the algorithmic workhorse for solving these optimization problems~\cite{boyd2004convex}. Due to their rapid convergence, rigorous theoretical foundations, and practical reliability, they remain the preferred choice over first-order methods—particularly when high-accuracy solutions are essential~\cite{nesterov1994interior,wright1997primal,nocedal2006numerical,potra2000interior}.

However, the robust performance of IPMs has historically required double-precision (f64) arithmetic to maintain the numerical stability of the linear systems solved at each iteration~\cite{wright2001finite,higham2002accuracy}. This requirement fundamentally clashes with modern machine learning hardware and software stacks, which prioritize single-precision (f32) and mixed-precision computation for throughput and memory efficiency~\cite{micikevicius2018mixed,jouppi2017datacenter}. While lower-precision formats accelerate parallelized pipelines, they severely amplify the numerical instabilities inherent in standard IPMs, rendering conventional solvers impractical for contemporary machine learning workflows~\cite{fasi2021numerical}.

Motivated by this critical bottleneck, we ask: \emph{how can we enable reliable, differentiable interior-point methods in single precision?} As a foundational step toward this broader goal, we focus on Quadratic Programs (QPs). QPs provide a canonical setting for studying these numerical issues, and while the ideas developed here extend to other convex cones, such extensions are beyond the scope of the present work.


\textbf{Contributions:} We introduce a differentiable interior-point method for QPs that (i) natively supports low-precision arithmetic and (ii) yields higher-accuracy gradients than standard IPM formulations. To realize this, we make the following contributions:
\begin{enumerate}[leftmargin=*]
    \item A novel primal-dual interior-point algorithm based on spectrally bounded linear systems that avoids the numerical instabilities affecting standard IPMs in low-precision arithmetic.
    \item A relaxation-based differentiation pipeline that leverages this algorithm to extract unique, smooth, and highly accurate derivatives of the solution with respect to the problem data.
    \item Extensive empirical validation on challenging benchmarks, demonstrating that our approach reliably solves and differentiates QPs in single precision where standard solvers strictly require double precision.
\end{enumerate}

\section{Preliminaries}
In this section, we introduce quadratic programming, the standard primal-dual interior-point framework, and its differentiability. This provides the foundations required to introduce our method.

\subsection{Quadratic Programming}

A quadratic program (QP) and its equivalent slack-variable formulation are given by
\begin{subequations}
\begin{center}
\begin{minipage}{0.42\textwidth}
\begin{equation}
\begin{aligned}
\underset{x}{\text{minimize}} \quad & \frac{1}{2} x^\top Q x + q^\top x \\
\text{subject to} \quad & Ax = b, \\
& Gx \leq h,
\end{aligned}
\end{equation}
\end{minipage}
\hfill $\Longleftrightarrow$ \hfill
\begin{minipage}{0.48\textwidth}
\begin{equation}
\begin{aligned}
\underset{x, s}{\text{minimize}} \quad & \frac{1}{2} x^\top Q x + q^\top x \\
\text{subject to} \quad & Ax = b, \\
& Gx + s = h, \\
& s \geq 0.
\end{aligned}
\end{equation}
\end{minipage}
\end{center}
\end{subequations}

Here, $x \in \mathbb{R}^n$ denotes the primal variable, $Q \in \mathbb{S}^n_+$ and $q \in \mathbb{R}^n$ define the quadratic objective, $A \in \mathbb{R}^{m \times n}$ and $b \in \mathbb{R}^m$ specify the equality constraints, $G \in \mathbb{R}^{p \times n}$ together with $h \in \mathbb{R}^p$ define the inequality constraints, and the slack variable $s \in \mathbb{R}^p$ converts inequalities into equalities.

Introducing Lagrange multipliers $y \in \mathbb{R}^m$ for the equality constraints and $z \in \mathbb{R}^p$ for the inequality constraints, the Lagrangian of the slack formulation is
\begin{equation}
\mathcal{L}(x, s, y, z) 
= \frac{1}{2} x^\top Q x + q^\top x 
+ y^\top (Ax - b) 
+ z^\top (Gx + s - h).
\end{equation}

The Karush--Kuhn--Tucker (KKT) conditions for optimality are given by
\begin{subequations}\label{eq:kkt_conditions}
\begin{alignat}{2}
& r_t \coloneqq \quad &Qx + q + G^\top z + A^\top y &= 0, \\
& r_c \coloneqq &z \odot s &= 0, \label{eq:kkt_conditions_complmentatirity}\\
& r_i \coloneqq &Gx + s - h &= 0, \\
& r_e \coloneqq &Ax -b &= 0, \\
&&s &\geq 0, \\
&&z &\geq 0,
\end{alignat}
\end{subequations}
where $\odot$ denotes elementwise multiplication and $r_t, r_c, r_e, r_i$ refer to the stationarity, complementarity, equality and inequality residuals. A primal--dual iterate $(x^*, s^*, y^*, z^*)$ is globally optimal if it satisfies~\eqref{eq:kkt_conditions}.

\subsection{Primal-Dual Interior-Point Methods}
To solve~\eqref{eq:kkt_conditions}, standard primal-dual IPMs perturb the complementarity condition in~\eqref{eq:kkt_conditions_complmentatirity} as
\begin{equation}
    r_{c,\kappa}\coloneqq z \odot s - \kappa,
\end{equation}
where $\kappa$ is the central path parameter. As $\kappa \rightarrow 0$, the solution of the perturbed system approaches that of~\eqref{eq:kkt_conditions}. For a fixed value of $\kappa$, Newton's method is applied to the modified KKT conditions, leading to the linear system
\begin{equation}\label{eq:explicit_uncondensed}
\begin{bmatrix}
Q & A^T & G^T & 0 \\
A & 0 & 0 & 0 \\
G & 0 & 0 & I \\
0 & 0 & D(s) & D(z)
\end{bmatrix}
\begin{bmatrix}
\Delta x \\
\Delta y \\
\Delta z \\
\Delta s
\end{bmatrix}
=
-
\begin{bmatrix}
r_t \\
r_e \\
r_i \\
r_{c,\kappa}
\end{bmatrix}.
\end{equation}
where $D(\cdot)$ denotes the diagonal matrix formed from a vector. The dominant computational cost at each iteration lies in solving this system, typically via a condensed formulation introduced later. 

Once the search directions are computed, a line search is performed to preserve the nonnegativity of $(s, z)$. The maximum step size $\alpha \leq 1$ such that $s + \alpha \Delta s \geq 0$ and $z + \alpha \Delta z \geq 0$ is given in closed form by
\begin{equation}\label{eq:linesearch}
\alpha = \min \left( 
1,\ 
\min_{i:\ \Delta s_i < 0} \frac{-s_i}{\Delta s_i},\ 
\min_{i:\ \Delta z_i < 0} \frac{-z_i}{\Delta z_i}
\right).
\end{equation}
For further details on primal-dual interior-point methods and their implementation, see~\cite{mattingley2012cvxgen}.

\subsection{Differentiable Interior-Point Quadratic Programming}\label{subsubsec:informative_gradients}
When differentiating numerical optimizers, it is computationally inefficient to propagate derivatives through the solver's internal iterations. An alternative approach is to solve the optimization problem nominally, and then compute the exact derivatives at the final solution using the implicit function theorem~\cite{krantz2012implicit} (see Appendix~\ref{appendix:implicit_function_theorem} for specific details on its applicability to interior-point QPs).
However, directly applying the implicit function theorem to differentiate a quadratic program can yield uninformative subgradients near active inequality constraints~\cite{agrawal2019differentiable}. To address this issue, derivatives are not evaluated at the strict optimum, where $\kappa^* \approx 0$, but instead at a relaxed point along the central path defined by $\kappa_{\mathrm{relax}} > \kappa^*$. This relaxation smooths the local geometry of the feasible set, producing gradients that remain well behaved, even near sharp constraint boundaries~\cite{tracy2024differentiability}. As a result, differentiating an interior-point QP naturally decomposes into three stages:
\begin{enumerate}
\item Solve the quadratic program to the desired convergence tolerance.
\item Relax the converged solution to a central path parameter $\kappa_{\mathrm{relax}}$.
\item Compute the exact derivatives at the relaxed point.
\end{enumerate}

The choice of $\kappa_{\mathrm{relax}}$ depends on the desired degree of smoothing. When only mild regularization is required, a small value is typically sufficient. In contrast, problems with strongly non-smooth geometry—such as sharp vertices or rapidly changing active sets—may benefit from a larger relaxation. Importantly, this approach decouples solution accuracy from gradient stability: the optimization problem can still be solved to high precision, while derivatives are computed from a smoother local representation. In modern automatic differentiation frameworks, these stages can be implemented naturally through customized forward and backward passes~\cite{jax2018github,NEURIPS2019_bdbca288}.

\section{Problem Statement}
The core computational bottleneck in both solving and differentiating interior-point quadratic programs lies in solving the linear system presented in~\eqref{eq:explicit_uncondensed}. To reduce the cost of matrix factorization, this full system is typically reduced to a more compact form:

\begin{equation}\label{eq:explicit_partially_condensed}
\begin{bmatrix}
Q & A^\top & G^\top \\
A & 0 & 0 \\
G & 0 & -D(s \oslash z)
\end{bmatrix}
\begin{bmatrix}
\Delta x \\
\Delta y \\
\Delta z
\end{bmatrix}
=
-
\begin{bmatrix}
r_t \\
r_e \\
r_i-D(z)^{-1} r_{c,\kappa}
\end{bmatrix},
\end{equation}
where $\oslash$ denotes the elementwise division operator.
This compact form exposes a fundamental numerical limitation of the standard primal-dual interior-point formulation: As the solver approaches the optimal solution, complementarity implies $s \odot z \rightarrow 0$. Consequently, the entries of the diagonal scaling matrix $D(s \oslash z)$ can become highly disparate. For inactive constraints, $z_i \rightarrow 0$ while $s_i$ remains positive, so $s_i/z_i$ diverges. For active constraints, $s_i \rightarrow 0$ while $z_i$ remains positive, so $s_i/z_i$ vanishes.

In double-precision arithmetic, this inherent ill-conditioning has been extensively studied and is generally considered manageable~\cite{wright1998ill, benzi2005numerical}, as it can be effectively handled using direct matrix factorizations. In single precision, however, this ill-conditioning introduces a severe risk of unstable or inaccurate factorizations, which either lead to poor step directions, unreliable gradients or, even worse, failures resulting from numerical overflows or underflows~\cite{abdelfattah2021survey,higham2022mixed}.
Furthermore, the additional linear solves required for differentiation increase exposure to the same numerical instabilities, compounding the risk of inaccurate gradients or solver failure. This numerical vulnerability makes robust single-precision differentiable interior-point solvers difficult, and has led many existing implementations to rely on double precision, preventing them from fully exploiting modern, low-precision machine learning accelerators~\cite{haidar2020mixed,agrawal2019differentiable,amos2017optnet, chen2024cuclarabel, tracy2024differentiability}.
To bridge this divide and enable robust differentiable single-precision interior-point solvers, the upcoming section introduces an alternative complementarity representation that effectively circumvents this structural ill-conditioning.

\section{Implicit Complementarity}~\label{sec:implicit_complementarity}
To satisfy the complementarity condition in~\eqref{eq:kkt_conditions_complmentatirity}, we adopt an approach similar to~\cite{permenter2023log, arrizabalaga2026implicit,bishop2026complementarity}, where instead of treating complementarity \emph{explicitly}, we satisfy it \emph{implicitly}, i.e., by construction.  We achieve this by introducing an auxiliary variable $v \in \mathbb{R}^p$ and a retraction map
\begin{equation}\label{eq:retraction_map}
    b_\kappa : \mathbb{R} \rightarrow \mathbb{R}_{+}\;\;\text{such that}\;\;b_\kappa(v) \cdot b_\kappa(-v) = \kappa\,,
\end{equation}
which relates the auxiliary variable by enforcing the following residuals to vanish: 
\begin{subequations}\label{eq:implicit_complementarity_residuals}
\begin{align}
    &r_z \coloneqq z - b_\kappa(v) = 0\,,\\
    &r_s \coloneqq s - b_\kappa(-v) = 0 \,.
\end{align}
\end{subequations}

\subsection{Choosing the retraction map}
From all the candidate retraction maps, we pick the \emph{softplus}: 
\begin{equation}\label{eq:softplus}
    b_\kappa(v) = \frac{v + \sqrt{v^2 + 4\kappa}}{2} \,.
\end{equation}
 This specific choice is motivated by the following two properties that its derivative holds
\begin{subequations}\label{eq:softplus_properties}
\begin{align}
&\partial_v b_\kappa(v) + \partial_v b_\kappa(-v) = 1\,,\\
&0 < \partial_v b_\kappa(v) \leq 1\,. \label{eq:softplus_properties_bounded_derivative}
\end{align}
These will play a fundamental role in the derivations of the underlying linear system that we will use in our approach.
\end{subequations}

\subsection{Newton's Method}
Because we operate on the manifold defined by $z \odot s = \kappa$, at each Newton step we target a specific central path parameter value, $\kappa_{\mathrm{target}}$. This introduces an additional residual:
\begin{equation}
    r_\kappa \coloneqq \kappa - \kappa_\text{target} = 0 \,.
\end{equation}
To accommodate the new variables, our method solves for the augmented direction $(\Delta x, \Delta y, \Delta z, \Delta s, \Delta v, \Delta \kappa)$. After the derivations in Appendix~\ref{subappendix:derivation_implicit_kkt}, the corresponding uncondensed KKT system is:

\begin{equation}\label{eq:implicit_kkt}
\begin{bmatrix}
    Q & A^\top & G^\top & 0 & 0 & 0 \\
    A & 0 & 0 & 0 & 0 & 0 \\
    G & 0 & 0 & I & 0 & 0 \\
    0 & 0 & I & 0 & -B_\kappa(v) & -c \\
    0 & 0 & 0 & I & B_\kappa(-v) & -c \\
    0 & 0 & 0 & 0 & 0 & 1
\end{bmatrix}
\begin{bmatrix}
    \Delta x \\ \Delta y \\ \Delta z \\ \Delta s \\ \Delta v \\ \Delta \kappa
\end{bmatrix}
=
-\begin{bmatrix}
    r_t \\ r_e \\ r_i \\ r_z \\ r_s \\ r_\kappa
\end{bmatrix}.
\end{equation}
where $B_\kappa(v) = \operatorname{Diag}\!\left(\partial_v b_\kappa(v)\right)$ and $c = \partial_\kappa b_\kappa(v)$. By exploiting the aforementioned properties of the retraction map in~\eqref{eq:softplus_properties}, we can reduce it into the following form\footnote{For alternative partially condensed forms see Appendix~\ref{subappendix:partially_condensed_ikkt}.}:
\begin{equation} \label{eq:implicit_partially_condensed}
\begin{bmatrix}
    Q-G^\top G & G^\top & A^\top \\
    G & -B_\kappa(-v) & 0 \\
    A & 0 & 0
\end{bmatrix}
\begin{bmatrix}
    \Delta x \\ \Delta v \\ \Delta y
\end{bmatrix}
=
-\begin{bmatrix}
    r_t - G^\top (r_{i}+r_z-r_s) \\
    r_{i}-r_s-c \, r_\kappa \\
    r_{e}
\end{bmatrix}.
\end{equation} 

The matrix of this reduced system exhibits two crucial properties: (i) it is symmetric, allowing for efficient $LDL^\top$ factorization~\cite{golub1996matrix}, and (ii) all of its block components remain constant across iterations, with the sole exception of $B_{\kappa}(-v)$. Notably, this middle block is diagonal and strictly bounded between $0$ and $1$ due to the derivative range of the retraction map~\eqref{eq:softplus_properties_bounded_derivative}. These properties yield two primary advantages:
First, the combination of symmetry and a predominantly static structure---where the only varying component is bounded and diagonal---enables a wide array of efficient linear solver strategies~\cite{arrizabalaga2026implicit}. Second, the eigenvalues of the linear system remain strictly bounded by the derivative range of the retraction map. This formulation circumvents the eigenvalue divergence that typically destabilizes standard interior-point systems in \eqref{eq:explicit_partially_condensed} as the iterates approach the optimal solution. Consequently, the linear system in \eqref{eq:implicit_partially_condensed} is structurally robust against the numerical instabilities that often hinder standard interior-point methods in low-precision arithmetic.
\begin{algorithm}[b]
\caption{Solving a Quadratic Program with Implicit  Complementarity}
\label{alg:solve_qp_implicit}
\small
\begin{algorithmic}[1]
\Function{solve\_qp}{$Q,q,A,b,G,h;\sigma,\mathrm{tol},\mathrm{max\_iter}$}
    \State $(x,s,z,y) \gets \Call{initialize}{Q,q,A,b,G,h}\quad \triangleright \text{same as in~\cite{vandenberghe2010cvxopt,mattingley2012cvxgen}}$

    \For{$i \gets 1:\mathrm{max\_iter}$}

        \Statex
        \State $\text{/* manifold coordinates */}$
        \State $v \gets z-s$
        \State $\kappa \gets \left(s^\top z\right)/m$

        \Statex
        \State $\text{/* residuals and check convergence */}$
        \State $r \gets \Call{evaluate\_residuals}{Q,q,A,b,G,h,x,y,z,s}\quad \triangleright \text{eqs.~\eqref{eq:kkt_conditions} and~\eqref{eq:implicit_complementarity_residuals}}$

        \State \textbf{if} $\|r\|_{\infty} < \mathrm{tol}$ \textbf{then} \Return $x,y,z,s$

        \Statex
        \State $\text{/* newton step */}$
        \State $\tilde{K}\gets\Call{precompute\_kkt\_factors}{Q,A,G,v,\kappa}\quad \triangleright \text{factorizations are cacheable}$
        \State $\kappa_{\mathrm{target}} \gets \sigma\kappa$
        \State $(\Delta x,\Delta y, \Delta z, \Delta s, \Delta v, \Delta\kappa)
        \gets
        \Call{solve\_kkt}{\tilde{K}, r,\kappa_{\mathrm{target}}}\quad \triangleright \text{eq.~\eqref{eq:implicit_kkt}}$

        \Statex
        \State $\text{/* line search and step on }(x,y,v,\kappa)\text{, then retract */}$
        \State $\alpha
        \gets
        \Call{linesearch}
        {s, z, \Delta s, \Delta z}\qquad \triangleright \text{eq.~\eqref{eq:linesearch}}$
        \State $x \gets x+\alpha\Delta x,\qquad
        y \gets y+\alpha\Delta y$
        \State $v \gets v+\alpha\Delta v,\qquad
        \kappa \gets \kappa+\alpha\Delta\kappa$
        \State $z, s \gets \Call{retraction\_map}{v,\kappa}\qquad \triangleright \text{eqs. ~\eqref{eq:implicit_complementarity_residuals} and~\eqref{eq:softplus}}$

    \EndFor
\EndFunction
\end{algorithmic}
\end{algorithm}

\section{Differentiable Implicit Primal-Dual Interior-Point Methods}
With a low-precision-compliant linear system established, we return to the differentiation procedure introduced in Section~\ref{subsubsec:informative_gradients}. As previously discussed, this three-step approach allows us to compute high-quality solutions to tight tolerances while evaluating smooth gradients at a relaxed point.

Here, we adapt this procedure to the implicit complementarity case and show how this sequence can be implemented as custom forward and backward passes compliant with modern automatic differentiation frameworks~\cite{jax2018github,NEURIPS2019_bdbca288}. 

\subsection{Solving the QP to high accuracy}

To solve the quadratic program to a tight tolerance, we integrate the implicit complementarity representation from the previous section into a primal-dual interior-point pipeline. An algorithmic overview is provided in Algorithm~\ref{alg:solve_qp_implicit}.
Two key aspects are worth highlighting: First, low-precision KKT solves must avoid the uncondensed form, relying exclusively on the partially condensed system in~\eqref{eq:implicit_partially_condensed} (see Appendix~\ref{subappendix:numerical_issues_ikkt}). Second, the target parameter $\kappa_\mathrm{target}$ is updated via a constant reduction factor $\sigma$. However, as in standard IPMs, this naturally can be extended to adaptive strategies, such as predictor-corrector schemes~\cite{mehrotra1992implementation}.

\subsection{Relaxing the QP solution}



The relaxation procedure is detailed in Algorithm~\ref{alg:relax_qp_implicit}. Its underlying mechanism closely mirrors that of Algorithm~\ref{alg:solve_qp_implicit}, with two primary distinctions: First, initialization is no longer required, as we warm-start the procedure directly from the optimal solution $(x^*,y^*,z^*,s^*)$ obtained via Algorithm~\ref{alg:solve_qp_implicit}. Second, rather than driving the central path parameter to zero ($\kappa \to 0$), we increase it toward a target relaxation value ($\kappa \to \kappa_\mathrm{relax}$). Because this relaxation essentially traces the central path in reverse, the matrix factorizations $\tilde{K}$ computed during the forward solve in Algorithm~\ref{alg:solve_qp_implicit} can be heavily reused.

\begin{algorithm}[h]
\caption{Relaxing an Implicit Quadratic Program to $\kappa_\mathrm{relax}$}
\label{alg:relax_qp_implicit}
\small
\begin{algorithmic}[1]
\Function{relax\_qp}{$Q,q,A,b,G,h,x,y,z,s,\kappa_\mathrm{relax};\tilde{K},\mathrm{tol},\mathrm{max\_iter}$}

    \For{$i \gets 1:\mathrm{max\_iter}$}

        \Statex
        \State $\text{/* manifold coordinates */}$
        \State $v \gets z-s$
        \State $\kappa \gets \left(s^\top z\right)/m$

        \Statex
        \State $\text{/* residuals and check convergence */}$
        \State $r \gets \Call{evaluate\_residuals}{Q,q,A,b,G,h,x,y,z,s}\quad \triangleright \text{eqs.~\eqref{eq:kkt_conditions} and~\eqref{eq:implicit_complementarity_residuals}}$

        \State \textbf{if} $\|r\|_{\infty} < \mathrm{tol}$ \textbf{then} \Return $x,y,z,s$

        \Statex
        \State $\text{/* newton step */}$
        \State \textbf{if} $\tilde{K}$ not cached $ \textbf{then}\; \tilde{K}\gets\Call{precompute\_kkt\_factors}{Q,A,G,v,\kappa}$
        \State $(\Delta x,\Delta y, \Delta z, \Delta s, \Delta v, \Delta\kappa)
        \gets
        \Call{solve\_kkt}{\tilde{K}, r, \kappa_{\mathrm{relax}}}\quad \triangleright \text{eq.~\eqref{eq:implicit_kkt}}$

        \Statex
        \State $\text{/* line search and step on }(x,y,v,\kappa)\text{, then retract */}$
        \State $\alpha
        \gets
        \Call{linesearch}
        {s, z, \Delta s, \Delta z}\qquad \triangleright \text{eq.~\eqref{eq:linesearch}}$
        \State $x \gets x+\alpha\Delta x,\qquad
        y \gets y+\alpha\Delta y$
        \State $v \gets v+\alpha\Delta v,\qquad
        \kappa \gets \kappa+\alpha\Delta\kappa$
        \State $z, s \gets \Call{retraction\_map}{v,\kappa}\qquad \triangleright \text{eqs.~\eqref{eq:implicit_complementarity_residuals} and~\eqref{eq:softplus}}$

    \EndFor
\EndFunction
\end{algorithmic}
\end{algorithm}

\subsection{Computing Gradients of the QP}
We compute the gradients at the relaxed point in a similar fashion to~\cite{amos2017optnet, tracy2024differentiability}. Instead of naively constructing the full Jacobian matrices—which can easily become prohibitively large—we apply the implicit function theorem to efficiently compute only the specific vector-Jacobian products required by the backward pass. The resulting procedure is detailed in Algorithm~\ref{alg:qp_gradients_implicit}.



Two key aspects of this procedure are worth noting: First, because the gradients are evaluated at $\kappa_\mathrm{relax}$, we can directly reuse the matrix factorizations computed in Algorithm~\ref{alg:relax_qp_implicit}. Second, the implicit reparameterization naturally absorbs the corresponding Jacobian factors, eliminating the need for the explicit multiplier scaling and unscaling steps required in standard primal-dual methods.

\begin{algorithm}[t]
\caption{Computing Gradients Through an Implicit Relaxed QP}
\label{alg:qp_gradients_implicit}
\small
\begin{algorithmic}[1]
\Function{compute\_qp\_gradients}{$Q,q,A,b,G,h,x,y,z,s,\kappa_\mathrm{relax}, \nabla_x \ell; \tilde{K}$}
    \State \text{/* compute differentials */}
    \State $r\gets (-\nabla_x \ell,0,0,0,0)$
    \State $(dx,dy,dz,ds,\_,\_)
    \gets
    \Call{solve\_kkt}
    {\tilde{K},r,\kappa_\mathrm{relax}}\qquad \triangleright \tilde{K} \text{cached in Alg.~\ref{alg:relax_qp_implicit}}$

    \Statex
    \State \text{/* parameter gradients */}
    \State $\nabla_Q \ell
    \gets
    \dfrac{(dx) x^\top + x (dx)^\top}{2}$

    \State $\nabla_q \ell \gets dx$

    \State $\nabla_A \ell
    \gets
    (dy) x^\top + y (dx)^\top$

    \State $\nabla_b \ell \gets -(dy)$

    \State $\nabla_G \ell
    \gets
    (dz) x^\top + z (dx)^\top$

    \State $\nabla_h \ell \gets -dz$

    \State \Return $\nabla_Q \ell,\nabla_q \ell,\nabla_A \ell,\nabla_b \ell,\nabla_G \ell,\nabla_h \ell$

\EndFunction
\end{algorithmic}
\end{algorithm}



\section{Numerical Experiments}\label{sec:numerical_experiments}
We evaluate our approach across three distinct benchmarks, structuring this section accordingly. In each experiment, we compare standard IPMs—where complementarity is handled \emph{explicitly} (visualized in blue)—against our \emph{implicit} formulation (visualized in orange). For clarity, we denote single-precision and double-precision arithmetic as f32 and f64, respectively.

To ensure a controlled and equitable comparison, we use the JAX-based framework of~\cite{tracy2023differentiable} as a baseline for state-of-the-art differentiable IPMs~\cite{amos2017optnet,chen2024cuclarabel,moreau2026}. We extend this codebase to incorporate our implicit method, allowing for an interchangeable evaluation of the two paradigms within the same underlying architecture. Experiments were performed on an NVIDIA GeForce RTX 3080 (12 GB). 

\subsection{Gradient Fidelity and Robustness in Single-Precision}
\begin{figure}[b]
  \centering
  \includegraphics[width=\linewidth]{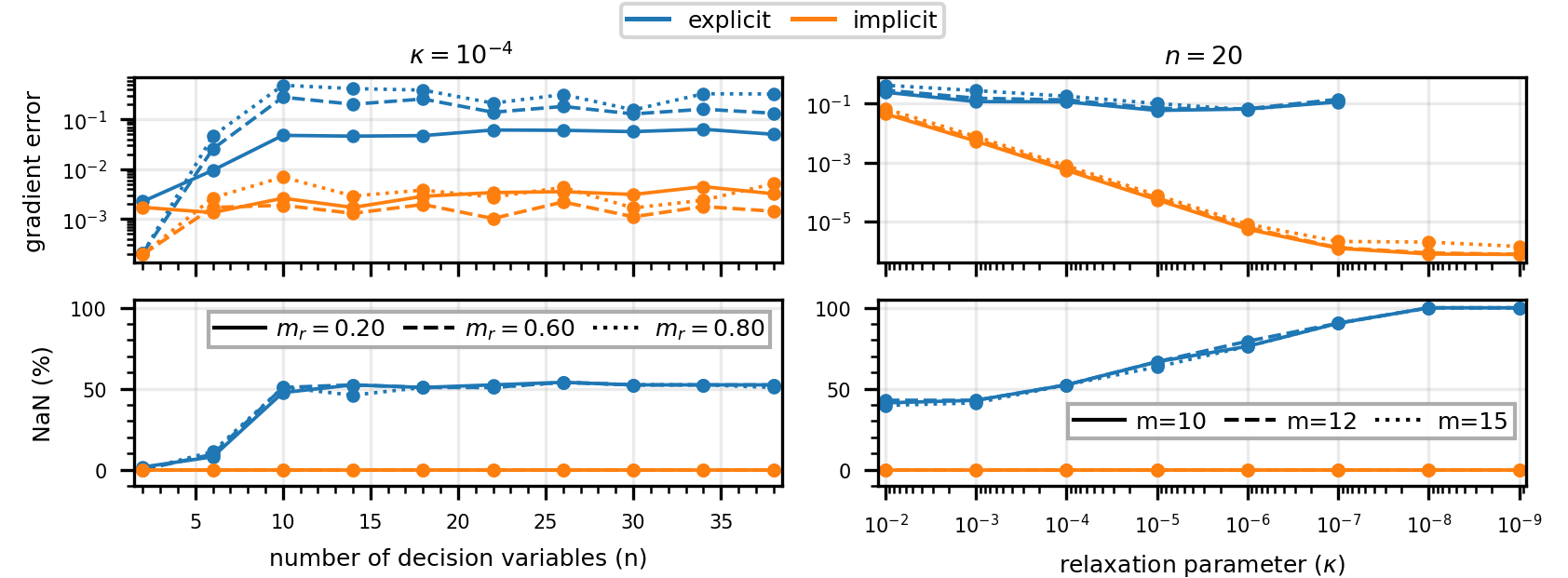}
  \caption{Implicit complementarity enables accurate and robust single-precision gradients. \textit{Left}: As problem dimension and active-set size increase, the implicit method maintains low gradient error and remains fully robust, whereas the explicit formulation fails 50\% of the time. Here, $m_r$ denotes the ratio of active inequality constraints $m$ to decision variables $n$. \textit{Right}: As $\kappa_{\mathrm{relax}}$ decreases, the implicit method continues to track the relaxation-limited gradient accuracy, whereas the explicit f32 formulation plateaus in accuracy and exhibits an escalating failure rate. Reported metrics are medians over three random seeds. }\label{fig:ablation} %
\end{figure}
In this experiment, we investigate the gradient quality and robustness of the implicit and explicit paradigms in single-precision (f32) arithmetic. Using the Euclidean projection onto a polytope against an exact active-set reference, we assess whether the projection map's derivative remains accurate and numerically stable as the solver relaxation parameter $\kappa_\mathrm{relax}$ varies, the active set grows, and problem dimensionality increases\footnote{See Appendix~\ref{appendix:ablation_study} for further details on the problem formulation and experiment design.}.

Results are depicted in Fig.~\ref{fig:ablation}. Sweeping over problem dimensions $n$ and constraints $m$ (first column), the implicit method's gradient error stays at the theoretical minimum dictated by the relaxation floor ($\kappa_\text{relax}=10^{-4}$). Conversely, the explicit method's accuracy degrades significantly with scale: for $n > 5$, it frequently triggers multiple active sets, causing severe numerical instability where 50\% of computations return NaNs (see Appendix~\ref{subappendix:nan_tracing} for a detailed trace). The implicit method, in contrast, maintains a 0\% failure rate and demonstrates total robustness to problem scale.

Evaluating the impact of the relaxation parameter $\kappa_\text{relax}$ at $n=20$ (second column) reveals that the explicit method's error stagnates near $10^{-1}$, while the implicit method faithfully tracks the $\kappa_\text{relax}$-allowed minimum until plateauing at $10^{-6}$. Moreover, as relaxation tightens, the explicit method becomes increasingly prone to failure, with NaN detection scaling from 50\% to 100\%.
These findings demonstrate the implicit approach's superior accuracy and reliability in single-precision settings, particularly when high-dimensional constraints or tight relaxations are required.

\subsection{Bilevel Optimization}
\begin{figure}
  \centering
  \includegraphics[width=\linewidth]{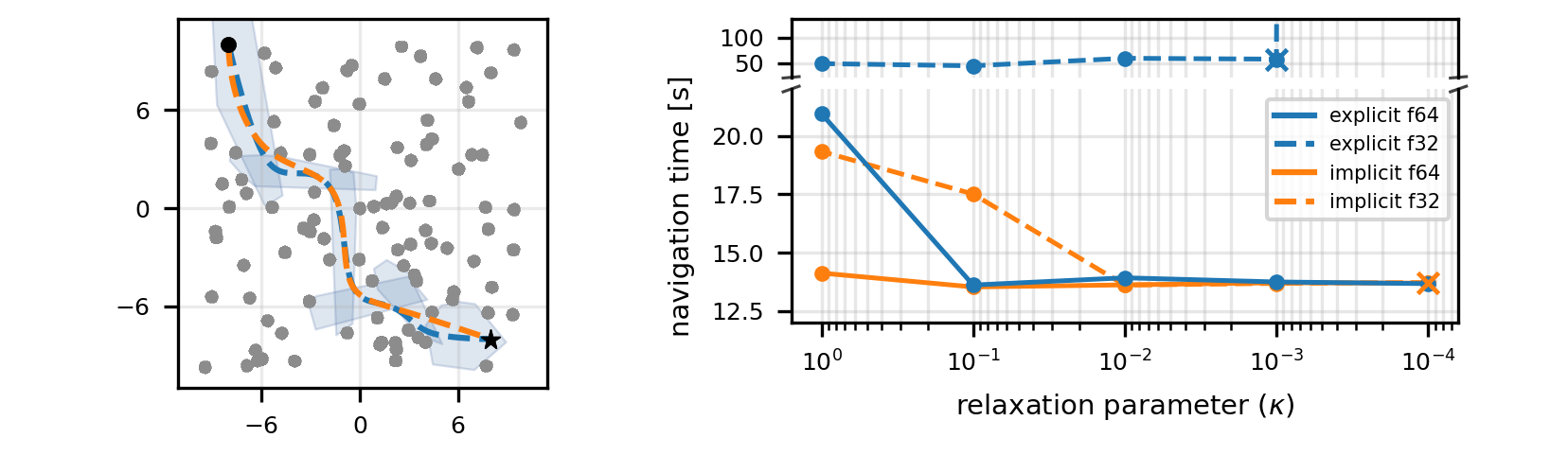}
    \caption{Bilevel trajectory optimization. Implicit f32 remains stable at tight relaxations and matches f64 performance, whereas explicit f32 becomes unreliable. \textit{Left:} Trajectories in forest environment. \textit{Right:} Navigation time for different relaxations. Crosses refer to the trajectories plotted in the left.}
  \label{fig:forest}
\end{figure}
Following our synthetic experiments, we evaluate our approach on a real-world bilevel trajectory optimization problem from robotics. The inner QP smooths the trajectory while enforcing kinematic and collision constraints, whereas the outer loop optimizes navigation time~\cite{mellinger2011minimum,liu2017planning,tordesillas2021faster}. Modeling the trajectory as a composite Bezier curve allows geometric safety to be enforced via linear constraints on control points, yielding a QP with coefficients that depend smoothly on segment durations (see Appendix~\ref{appendix:bilevel_opt}). For a fixed sequence of durations $T = (T_1,\dots,T_K)$, the inner QP minimizes a weighted smoothness objective subject to safety, boundary, and derivative constraints. The outer level then updates $T$ to minimize the inner objective alongside a linear time penalty. Because this outer objective depends directly on the inner QP's solution, optimizing $T$ requires differentiating through the solver.

We optimize the outer loop via L-BFGS~\cite{byrd1995limited,nocedal2006numerical}, comparing the usefulness of gradients provided by the explicit and implicit paradigms. We sweep across multiple relaxation parameters in single and double precision for two environments: a forest and an office. Since both scenarios yield similar behavior, we present the forest environment here (Fig.~\ref{fig:forest}) and defer the office to the appendix (Fig.~\ref{fig:office}).
These sweeps reveal two key insights: First, while excessively loose relaxations prevent optimal convergence, tightening the parameter beyond a certain threshold offers diminishing returns. Second, although both methods perform equivalently in double precision, the explicit method in single precision produces unusable gradients or fails entirely at tighter relaxations ($\kappa_\mathrm{relax}=10^{-4}$). In contrast, the implicit method in single precision remains robust, performing on par with its double-precision counterpart.
Aligning with our synthetic case study, these results confirm that the implicit method consistently delivers high-quality, reliable gradients across both single and double precision.

\subsection{Learning a Multi-Agent Safety Filter}
\begin{figure}
  \centering
  \includegraphics[width=\linewidth]{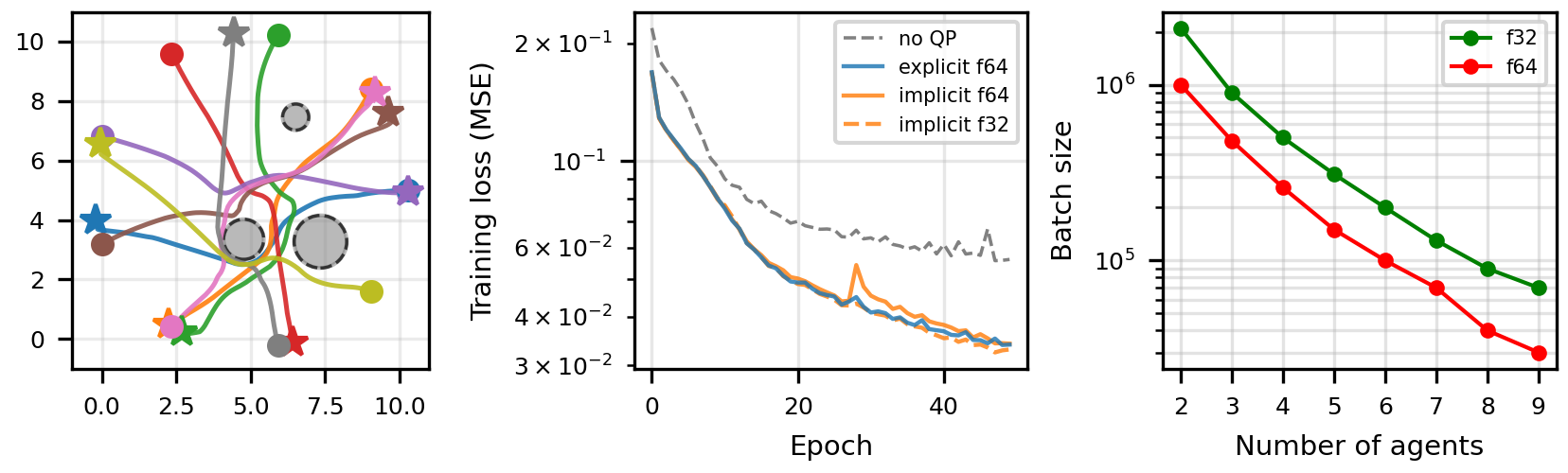}
  \caption{End-to-end learning of a multi-agent safety filter. Implicit f32 matches f64 training stability (unlike explicit f32, which fails entirely) and allows for double the maximum batch size.}
  \label{fig:safety_filter}
\end{figure}
While the bilevel trajectory optimization serves as a compelling application, one might argue it represents a specialized use case. In standard deep learning settings, stochastic optimizers such as SGD~\cite{robbins1951stochastic} and Adam~\cite{kingma2014adam} are often robust to noisy gradients, reducing the need for highly accurate derivatives~\cite{bottou2018optimization}. This raises a natural question: are explicit single-precision (f32) gradients truly insufficient, and does the improved gradient quality and robustness of our implicit approach offer a tangible advantage in end-to-end learning pipelines?

To investigate this question, we consider a centralized multi-agent safety-filter learning task. Here, a neural network predicts a nominal control law alongside state-dependent control barrier function (CBF) gains, and a differentiable quadratic programming (QP) layer projects this proposal onto a safe set (see Appendix~\ref{appendix:safety_filter} for formulation and implementation details)~\cite{xiao2023barriernet}. This environment is an ideal testbed because it directly addresses the aforementioned questions while allowing us to evaluate the high-dimensional regime explored in our synthetic case study; specifically, increasing the number of agents scales the QP's decision variables linearly and its inequality constraints quadratically.

We evaluate the learning process scaling from 1 to 9 agents in an environment with 3 obstacles. The results are summarized in Fig.~\ref{fig:safety_filter}. The first column illustrates an exemplary scenario with 9 agents, while the second displays the corresponding training loss curves. As a baseline, we also include a vanilla MLP lacking a QP safety filter. Notably, the implicit f32, implicit f64, and explicit f64 formulations exhibit near-identical, stable performance. Crucially, however, the explicit f32 method is omitted from these plots because it completely fails (NaNs) for any scenario involving more than a single agent. This contrast demonstrates how our implicit approach successfully unlocks the use of single-precision arithmetic in end-to-end learning with differentiable QPs—a regime where state-of-the-art methods strictly mandate double precision.

This capability offers a significant practical advantage, as modern machine learning hardware is natively optimized for single-precision operations. As illustrated in the third column of Fig.~\ref{fig:safety_filter}, the implicit f32 method accommodates GPU batch sizes exactly twice as large as its f64 counterpart, definitively demonstrating superior memory efficiency and scalability.


\section{Conclusions}
In this paper, we introduced a novel differentiable interior-point method that resolves the numerical instabilities of solving and differentiating quadratic programs in single-precision arithmetic. By replacing standard explicit complementarity with an implicit representation, our approach avoids the severe ill-conditioning that has historically confined differentiable optimization to double precision. Extensive evaluations—spanning synthetic benchmarks, bilevel trajectory optimization, and multi-agent safety filter learning—demonstrate that our implicit method matches the gradient accuracy and robustness of double-precision solvers. Crucially, by unlocking reliable low-precision operations, our formulation fully leverages the memory and throughput advantages of modern hardware, establishing a scalable foundation for embedding constrained optimization into machine learning pipelines.

\newpage
\bibliographystyle{unsrt}
\bibliography{references}

\newpage
\appendix
\section{Implicit Function Theorem in Primal-Dual Interior-Point QPs}\label{appendix:implicit_function_theorem}
Consider variables $w \in \mathbb{R}^a$ and parameters $\theta \in \mathbb{R}^b$, and suppose that an equilibrium point $w^*$ is defined implicitly by
\begin{equation}
    r(w^*, \theta) = 0.
\end{equation}
Linearizing the residual $r$ around $w^*$ gives
\begin{equation}
    \frac{\partial r}{\partial w} \delta w
    +
    \frac{\partial r}{\partial \theta} \delta \theta
    =
    0.
\end{equation}
Solving this equation for the parameter derivatives yields:
\begin{equation}\label{eq:implicit_function_theorem}
    \frac{\partial w}{\partial \theta}
    =
    -
    \left(
    \frac{\partial r}{\partial w}
    \right)^{-1}
    \frac{\partial r}{\partial \theta}.
\end{equation}

Applying this to the present formulation, the KKT conditions in~\eqref{eq:kkt_conditions} are treated as the residual $r$, the primal-dual variables $(x, s, z, y)$ as $w$, and the problem data as the parameters $\theta$. The implicit function theorem can therefore be used to differentiate the optimizer with respect to the problem data. Importantly, the linear system associated with~\eqref{eq:implicit_function_theorem} has the same structure as the Newton system in~\eqref{eq:explicit_uncondensed}. As a result, these derivatives can be computed efficiently by reusing the same linear-algebra framework, and potentially the same matrix factorizations, already available inside the interior-point solver.

\section{The KKT system of a QP with implicit complementarity}

\subsection{Derivation of the uncondensed form}\label{subappendix:derivation_implicit_kkt}
We start by linearizing the implicit complementarity residuals from~\eqref{eq:implicit_complementarity_residuals} with respect to $v$ and $\kappa$:
\begin{subequations}
\begin{align}
    r_z^+ &\approx r_z + \Delta z - B_\kappa(v) \, \Delta v - c \, \Delta \kappa\,, \\
    r_s^+ &\approx r_s + \Delta s + B_\kappa(-v) \Delta v - c \Delta \kappa\,, \\
    r_\kappa^+ &\approx r_\kappa + \Delta \kappa\,,
\end{align}
\end{subequations}
where $B_\kappa(v) = \operatorname{Diag}\!\left(\partial_v b_\kappa(v)\right)$ and $c = \partial_\kappa b_\kappa(v)$. Combining these residuals with the ones in~\eqref{eq:kkt_conditions} yields the following augmented linear system:
\begin{equation}\label{eq:implicit_kkt_appendix}
\begin{bmatrix}
    Q & A^\top & G^\top & 0 & 0 & 0 \\
    A & 0 & 0 & 0 & 0 & 0 \\
    G & 0 & 0 & I & 0 & 0 \\
    0 & 0 & I & 0 & -B_\kappa(v) & -c \\
    0 & 0 & 0 & I & B_\kappa(-v) & -c \\
    0 & 0 & 0 & 0 & 0 & 1
\end{bmatrix}
\begin{bmatrix}
    \Delta x \\ \Delta y \\ \Delta z \\ \Delta s \\ \Delta v \\ \Delta \kappa
\end{bmatrix}
=
-\begin{bmatrix}
    r_t \\ r_e \\ r_i \\ r_z \\ r_s \\ r_\kappa
\end{bmatrix}.
\end{equation}

\subsection{Numerical issues of the uncondensed form in low-precision arithmetic}
\label{subappendix:numerical_issues_ikkt}
The uncondensed system in eq.~\eqref{eq:implicit_kkt_appendix} is problematic in low-precision arithmetic. This is because the auxiliary increment $\Delta v$ influences the Newton step primarily through its coupling with $\Delta z$ via $B_\kappa(v)$. Because $\Delta z$ enters the stationarity equation directly, this indirect dependency becomes a critical bottleneck. In particular, as $\kappa$ vanishes and a constraint becomes inactive, the diagonal entries of $B_\kappa(v)$ approach zero. Under low-precision regimes, this causes the signal from $\Delta v$ to be effectively lost, leading to a breakdown in the coupling required for accurate stationarity updates.

\subsection{Partially condensed forms}
\label{subappendix:partially_condensed_ikkt}

The uncondensed system can be reduced into two distinct partially condensed forms. The first and most intuitive reduction follows the standard IPM structure:
\begin{equation}
\begin{bmatrix} 
    Q & G^\top B_\kappa(v) & A^\top \\
    G & -B_\kappa(-v) & 0 \\ 
    A & 0 & 0 
\end{bmatrix} 
\begin{bmatrix} 
    \Delta x \\ \Delta v \\ \Delta y 
\end{bmatrix} 
=
- \begin{bmatrix} 
    r_t - G^\top(r_z + c\,r_\kappa) \\ 
    r_{i}-r_s-c\,r_\kappa \\ 
    r_{e} 
\end{bmatrix}.
\end{equation}
While this system closely resembles the traditional formulation in eq.~\eqref{eq:explicit_partially_condensed}, it remains explicitly dependent on $B_\kappa(v)$. Consequently, it inherits the numerical instabilities detailed in Appendix~\ref{subappendix:numerical_issues_ikkt}, making it unreliable for low-precision computations.

To address these limitations, we propose the second candidate introduced in eq.~\eqref{eq:implicit_partially_condensed}:
\begin{equation}
	\begin{bmatrix}
    	Q-G^\top G & G^\top & A^\top \\
    	G & -B_\kappa(-v) & 0 \\
    	A & 0 & 0
	\end{bmatrix}
	\begin{bmatrix}
    	\Delta x \\ \Delta v \\ \Delta y
	\end{bmatrix}
	= - \begin{bmatrix}
    	r_t - G^\top(r_{i}+r_z-r_s) \\
    	r_{i}-r_s-c\,r_\kappa \\
    	r_{e}
	\end{bmatrix}.
\end{equation}
Although this approach densifies the Hessian, it offers the significant advantages discussed in Section~\ref{sec:implicit_complementarity}: the system is symmetric, predominantly static, and spectrally bounded. Most importantly, by eliminating the reliance on $B_\kappa(v)$, it bypasses the numerical issues of the uncondensed form, establishing itself as the most robust candidate for low-precision arithmetic.

\section{Implementation of the derivative of the softplus retraction map}\label{appendix:di2pd_tricks}
To avoid catastrophic cancellation (when $v<0$ and $v^2>>\kappa$) in the retraction map $b_\kappa(v)$ and its derivative $\partial_v b_\kappa(v$) , the software implementation of the softplus retraction map and its derivative needs to be carefully conducted:
\begin{subequations}
    \begin{align}
b_\kappa(v) &= 
\begin{cases} 
\frac{v + \sqrt{v^2 + 4\kappa}}{2} & \text{if } v \ge 0 \\ 
\frac{2\kappa}{\sqrt{v^2 + 4\kappa} - v} & \text{if } v < 0 
\end{cases} \\[10pt]
\partial_v b_\kappa(v) &= 
\begin{cases} 
\frac{1}{2} \left( 1 + \frac{v}{\sqrt{v^2 + 4\kappa}} \right) & \text{if } v \ge 0 \\ 
\frac{2\kappa}{v^2 + 4\kappa - v\sqrt{v^2 + 4\kappa}} & \text{if } v < 0 
\end{cases}
\end{align}
\end{subequations}

\section{Gradient Fidelity and Robustness in Single-Precision}\label{appendix:ablation_study}
\subsection{Problem Formulation}
Let
\begin{equation}
  \mathcal{P} = \{ y \in \mathbb{R}^n : Gy \leq h \}
\end{equation}
be a polytope defined by $p$ linear inequalities, where $G \in \mathbb{R}^{p \times n}$ and $h \in \mathbb{R}^{p}$. For a query point $x \in \mathbb{R}^n$, the experiment studies the Euclidean projection
\begin{equation}
  \Pi_{\mathcal{P}}(x)
  \;=\;
  \arg\min_{y \in \mathbb{R}^n}
  \frac{1}{2}\|y - x\|_2^2
  \quad
  \text{s.t.}
  \quad
  Gy \leq h.
  \label{eq:projection}
\end{equation}
This problem is encoded into the QPs with the following problem data
\begin{equation}
  Q = I_n,\qquad c = -x,\qquad A_{\mathrm{eq}} \in \mathbb{R}^{0 \times n},\qquad b_{\mathrm{eq}} \in \mathbb{R}^{0},
\end{equation}
so the objective solved in standard QP form is
\begin{equation}
  \frac{1}{2} y^\top y - x^\top y,
\end{equation}
which differs from \eqref{eq:projection} only by the constant $-\frac{1}{2}\|x\|_2^2$.

The quantity of interest is the gradient of the projection map with respect to the query point. Rather than comparing full Jacobian matrices, we benchmark the output with a fixed probe direction $v \in \mathbb{R}^n$ and differentiate
\begin{equation}
  \phi(x) = v^\top \Pi_{\mathcal{P}}(x).
\end{equation}
For each sampled instance, we compute
\begin{equation}
  g(x) = \nabla_x \phi(x)
\end{equation}
through automatic differentiation of the QP solver.

\subsection{Exact Hard-Projection Reference}
The reference gradient is constructed analytically from the active set. The generator builds instances so that the true projection point is a prescribed point $y_0$ on a face of the polytope with exactly $m$ active constraints. Let $G_A \in \mathbb{R}^{m \times n}$ denote the corresponding active rows. 

Under the fixed-active-set model, the local projection map is the orthogonal projection onto the affine subspace
\begin{equation}
  \{ y : G_A y = h_A \},
\end{equation}
whose Jacobian is defined as
\begin{equation}
  J_{\mathrm{hard}} = I_n - G_A^\top (G_A G_A^\top)^{-1} G_A.
  \label{eq:hard-jacobian}
\end{equation}
Equation~\eqref{eq:hard-jacobian} represents the orthogonal projector onto the tangent space of the active face. Since $J_{\mathrm{hard}}$ is symmetric, the exact gradient of the scalarized objective $\phi(x) = v^\top \Pi_{\mathcal{P}}(x)$ is given by
\begin{equation}
  g_{\mathrm{hard}} = J_{\mathrm{hard}} v.
\end{equation}

\subsection{Experiment design}\label{subapp:ablation_experimental_design}


To study how geometry influences gradient quality, we vary the location of the query point relative to the polytope boundary. More precisely, we control how far the point lies outside the polytope while keeping the projected point on the same active face. This isolates the effect of moving away from the boundary along the outward normal directions.

For each instance, we first sample a boundary point $y_0$ and a set of $m$ active constraints, collected in the matrix $G_A$. We then define the query point as
\begin{equation}
  x = y_0 + G_A^\top \lambda,
  \qquad
  \lambda = d\,\xi,
\end{equation}
where the components of $\xi$ are sampled independently as $\xi_i \sim \mathrm{Uniform}(0.5, 1.5)$ and $d \in \mathrm{logspace}(-2, 2, 21)$. In this construction, the scalar $d$ controls the overall displacement away from the boundary point $y_0$. As $d$ increases, the query point moves farther outward from the active face, and both the constraint violation and the distance $\|x-y_0\|_2$ increase approximately linearly.

Across all experiments we maintain the following fixed:
\begin{align}
  &\mathrm{solver\_max\_iter} = 100, \\
  &\text{seed} \in \{0,1,2\}.
\end{align}

\textbf{Sweep over problem size}. We fix $\kappa_\mathrm{relax}=10^{-4}$ and sweep the following:
\begin{align}
  &n \in \{2, 6, 10, 14, 18, 22, 26, 30, 34, 38\}, \\
  &p = 1.25n , \\
  &m_r \in \{0.2, 0.6, 0.8\},\\
  &\mathrm{solver\_tol} = 10^{-4}.
\end{align}

\textbf{Sweep over relaxation parameter}. We fix $n=20,p=25,m=\{10,12,15\}$ and sweep over relaxation parameter as:
\begin{align}
  &\kappa_\mathrm{relax} \in \{10^{-2},10^{-3},10^{-4},10^{-5},10^{-6},10^{-7},10^{-8},10^{-9}\}\\
  &\mathrm{solver\_tol} = \min(\kappa_\mathrm{relax}, 10^{-4})
\end{align}

\subsection{NaN tracing in explicit f32}~\label{subappendix:nan_tracing}
For the $\kappa$ sweep in Fig.~\ref{fig:ablation}, we trace the exact operations where first non-finite intermediate appears inside f32 explicit operations.

Over the full sweep, first-hit failures are attributed only to the affine predictor step of the primal-dual interior-point method or to the subsequent relaxed continuation; all other traced groups remain inactive.

\begin{figure}
  \centering
  \includegraphics[width=\linewidth]{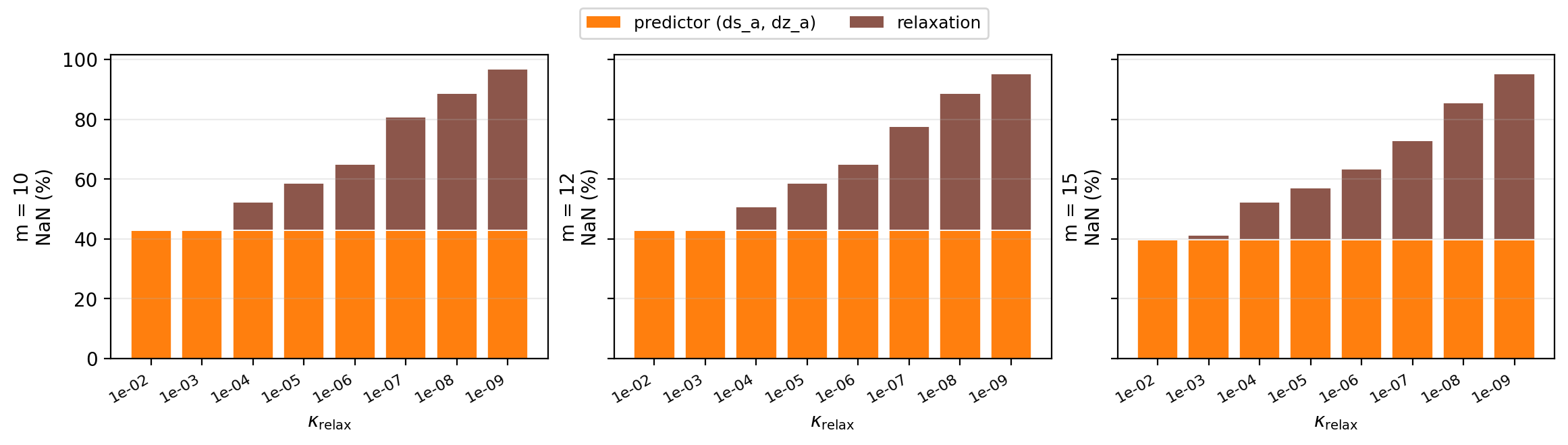}
  \caption{NaN tracing for operations in explicit f32.}~\label{fig:nan_tracing}
\end{figure}

In particular we trace the following operations in the codebase of~\cite{tracy2023differentiable}. First, in the forward primal-dual interior-point solve, we separately test the KKT scaling vector $z/s$, the affine predictor directions $(ds_a,dz_a)$, the centering quantities $(\sigma,\mu)$, the corrector directions $(dx,ds,dz,dy)$, and the line-search step $\alpha$. Then, we also check the relaxed solve and the custom backward solve. The results are given in the table~\ref{table:nan_tracing} and plotted in Fig.~\ref{fig:nan_tracing}.

\begin{table}[t]
  \centering
  \caption{Operation-level first-hit nan tracing. Each row summarizes $63$ instances. Column abbreviations are: Sc.\ (scaling), Pred.\ (predictor), Cent.\ (centering), Corr.\ (corrector), LS (line search), Rel.\ (relaxation), Bwd.\ (backward), and N/A (no attribution).}
  \label{table:nan_tracing}
  \scriptsize
  \resizebox{\textwidth}{!}{%
  \begin{tabular}{ccccccccccc}
    \toprule
    $m$ & $\kappa$ & total & Sc. & Pred. & Cent. & Corr. & LS & Rel. & Bwd. & N/A \\
    \midrule
    10 & $10^{-2}$ & 63 & 0 & 26 & 0 & 0 & 0 & 0  & 0 & 0 \\
    10 & $10^{-3}$ & 63 & 0 & 26 & 0 & 0 & 0 & 1  & 0 & 0 \\
    10 & $10^{-4}$ & 63 & 0 & 26 & 0 & 0 & 0 & 7  & 0 & 0 \\
    10 & $10^{-5}$ & 63 & 0 & 26 & 0 & 0 & 0 & 12 & 0 & 0 \\
    10 & $10^{-6}$ & 63 & 0 & 26 & 0 & 0 & 0 & 16 & 0 & 0 \\
    10 & $10^{-7}$ & 63 & 0 & 26 & 0 & 0 & 0 & 22 & 0 & 0 \\
    10 & $10^{-8}$ & 63 & 0 & 26 & 0 & 0 & 0 & 28 & 0 & 0 \\
    10 & $10^{-9}$ & 63 & 0 & 26 & 0 & 0 & 0 & 33 & 0 & 0 \\
    12 & $10^{-2}$ & 63 & 0 & 27 & 0 & 0 & 0 & 0  & 0 & 0 \\
    12 & $10^{-3}$ & 63 & 0 & 27 & 0 & 0 & 0 & 0  & 0 & 0 \\
    12 & $10^{-4}$ & 63 & 0 & 27 & 0 & 0 & 0 & 6  & 0 & 0 \\
    12 & $10^{-5}$ & 63 & 0 & 27 & 0 & 0 & 0 & 10 & 0 & 0 \\
    12 & $10^{-6}$ & 63 & 0 & 27 & 0 & 0 & 0 & 15 & 0 & 0 \\
    12 & $10^{-7}$ & 63 & 0 & 27 & 0 & 0 & 0 & 22 & 0 & 0 \\
    12 & $10^{-8}$ & 63 & 0 & 27 & 0 & 0 & 0 & 26 & 0 & 0 \\
    12 & $10^{-9}$ & 63 & 0 & 27 & 0 & 0 & 0 & 33 & 0 & 0 \\
    15 & $10^{-2}$ & 63 & 0 & 25 & 0 & 0 & 0 & 0  & 0 & 0 \\
    15 & $10^{-3}$ & 63 & 0 & 25 & 0 & 0 & 0 & 1  & 0 & 0 \\
    15 & $10^{-4}$ & 63 & 0 & 25 & 0 & 0 & 0 & 7  & 0 & 0 \\
    15 & $10^{-5}$ & 63 & 0 & 25 & 0 & 0 & 0 & 10 & 0 & 0 \\
    15 & $10^{-6}$ & 63 & 0 & 25 & 0 & 0 & 0 & 14 & 0 & 0 \\
    15 & $10^{-7}$ & 63 & 0 & 25 & 0 & 0 & 0 & 20 & 0 & 0 \\
    15 & $10^{-8}$ & 63 & 0 & 25 & 0 & 0 & 0 & 28 & 0 & 0 \\
    15 & $10^{-9}$ & 63 & 0 & 25 & 0 & 0 & 0 & 36 & 0 & 0 \\
    \bottomrule
  \end{tabular}}
\end{table}

\section{Bilevel Trajectory Optimization} 
\label{appendix:bilevel_opt}

\subsection{Problem Formulation}

\subsubsection{Bézier curve parameterization}
Let \(K\) denote the number of convex safe cells, and let each trajectory segment be a degree-\(n\) Bézier curve in \(\mathbb{R}^d\). For our sweep experiments, we use \(n=4\) and \(d=2\). Each segment \(k\) is parameterized over a normalized local coordinate \(u \in [0,1]\) as:
\begin{equation}
p_k(u) = \sum_{j=0}^{n} P_{k,j} b_j^{n}(u),
\end{equation}
where \(P_{k,j} \in \mathbb{R}^d\) are the control points and \(b_j^n\) are the Bernstein polynomials. The physical time duration assigned to segment \(k\) is \(T_k > 0\), yielding a total traversal time of:
\begin{equation}
T_{\mathrm{tot}} = \sum_{k=1}^{K} T_k.
\end{equation}
To compute trajectory derivatives efficiently, for each derivative order \(r\), we construct a differentiation matrix \(M_r\) such that the control points of the \(r\)-th derivative with respect to the normalized parameter \(u\) are simply \(M_r P_k\). When mapping back to physical time \(\tau\), the derivative must be scaled by \(T_k^{-r}\):
\begin{equation}
\frac{d^r p_k}{d\tau^r}(\tau) = \frac{1}{T_k^r} \frac{d^r p_k}{du^r}(u).
\end{equation}
This temporal scaling is the core source of the bilevel coupling in our formulation: both the quadratic cost matrix and several of the kinematic constraints in the inner problem depend non-linearly on the duration vector \(T\).

\subsubsection{Inner level: quadratic program over control points}
Given a fixed segment duration vector \(T\), we formulate the inner optimization problem by stacking all control points into a single decision vector \(x \in \mathbb{R}^{K(n+1)d}\). The resulting trajectory optimization problem is a purely quadratic program (QP):
\begin{equation}
\begin{aligned}
\phi(T) = \min_{x} \quad & \frac{1}{2} x^\top Q(T) x \\
\textrm{s.t.} \quad & A_{\mathrm{eq}}(T) x = b_{\mathrm{eq}}(T), \\
& Gx \le h(T).
\end{aligned}
\label{eq:inner_qp}
\end{equation}
Note that the linear objective term is zero in our implementation, so the objective is entirely quadratic. 

\paragraph{Objective Function.} The cost matrix \(Q(T)\) is block-diagonal across the trajectory segments and is designed to penalize the derivatives of the curve. For a given segment \(k\), derivative order \(r\), and penalty weight \(w_r\), the continuous-time cost is:
\begin{equation}
\frac{w_r}{T_k^{2r-1}} \int_0^1 \left\lVert \frac{d^r p_k}{du^r}(u) \right\rVert^2 du.
\end{equation}
Using Bernstein Gram matrices, this integral is evaluated exactly and assembled as a quadratic form over the control points. For our sweep experiments, we penalize velocity, acceleration, jerk, and snap using the weights:
\begin{equation}
(w_1, w_2, w_3, w_4) = (0.1, \, 0.1, \, 1.0, \, 0.1).
\end{equation}
This configuration heavily penalizes jerk while applying milder regularization to the other derivatives.

The equality and inequality constraints in \eqref{eq:inner_qp} enforce safety, continuity, and boundary conditions as follows:

\paragraph{Cell-membership constraints.} We restrict the trajectory to safe regions \footnote{We obtain the convex safe regions by using the python package \href{https://pypi.org/project/pydecomp/}{pydecomp}.} defined by convex polytopes\(\mathcal{C}_k = \{p \in \mathbb{R}^2 \mid A_k p \le b_k\}\). By constraining every Bézier control point of segment \(k\) to lie within \(\mathcal{C}_k\), the convex-hull property of Bézier curves guarantees that the entire segment remains safely within the cell.

\paragraph{Endpoint constraints.} The initial control point of the first segment is fixed to the designated start state, and the final control point of the last segment is fixed to the goal state.

\paragraph{Continuity constraints.} To ensure the physical trajectory is \(C^2\) continuous, adjacent segments must match up to the second derivative (\(c=2\)). For \(r \in \{0, 1, 2\}\), we enforce:
\begin{equation}
\frac{(M_r P_k)_{\mathrm{end}}}{T_k^r} = \frac{(M_r P_{k+1})_{\mathrm{start}}}{T_{k+1}^r}.
\end{equation}

\paragraph{Initial and terminal derivative constraints.} The trajectory must start and end at rest. We enforce zero velocity and acceleration at the boundaries for \(r \in \{1, 2\}\):
\begin{equation}
\frac{(M_r P_1)_{\mathrm{start}}}{T_1^r} = 0, \qquad \frac{(M_r P_K)_{\mathrm{end}}}{T_K^r} = 0.
\end{equation}

\paragraph{Derivative bounds.} We impose strict velocity and acceleration limits by bounding the derivative control polygons. Because the derivative of a Bézier curve is itself a Bézier curve, constraining its control points strictly bounds the derivative along the entire segment. For \(r \in \{1,2\}\), we enforce:
\begin{equation}
\ell_r \le \frac{M_r P_k}{T_k^r} \le u_r.
\end{equation}
In our sweep setup, the velocity limits are \(\ell_1 = (-10, -10)\) and \(u_1 = (10, 10)\), while the acceleration limits are \(\ell_2 = (-2, -2)\) and \(u_2 = (2, 2)\).
\begin{figure}[t]
  \centering
  \includegraphics[width=\linewidth]{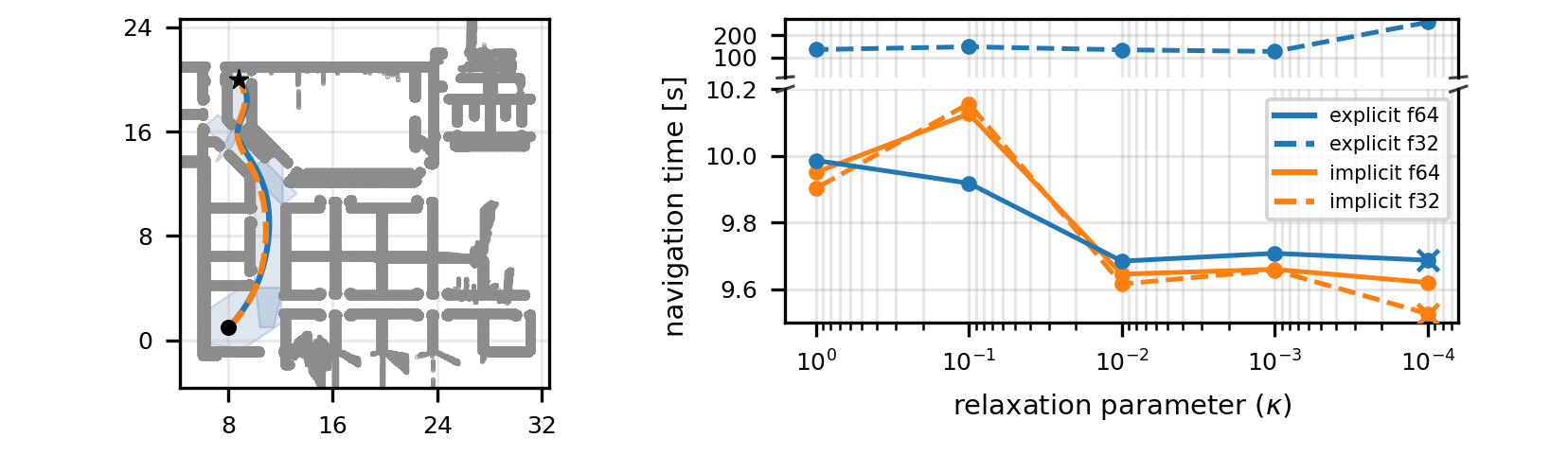}
  \caption{Bilevel optimization for fast (outer opt.) and safe (inner opt.) trajectories. \textit{Left:} Trajectories in office environment. \textit{Right:} Navigation time for trajectories computed with different relaxations. Crosses in the right column plot refer to the trajectories picked for plotting in the left column.}\label{fig:office}
\end{figure}
\subsubsection{Outer level: optimization over segment durations}
The outer optimization seeks to find the optimal segment durations by augmenting the inner smoothness cost \(\phi(T)\) with a linear time penalty:
\begin{equation}
\min_{T \in \mathbb{R}^K_{>0}} \mathcal{F}(T) := \phi(T) + \lambda \sum_{k=1}^{K} T_k \quad \text{s.t.} \quad T_k \ge T_{\mathrm{min}}.
\end{equation}
In our experiments, we set \(\lambda = 100\) and \(T_{\mathrm{min}} = 0.1\). Without this linear time penalty, the optimization would naturally favor arbitrarily large durations to minimize the smoothness cost. By penalizing time, the problem becomes genuinely bilevel: reducing the duration \(T_k\) decreases the linear penalty but simultaneously tightens the feasible set and alters the objective curvature of the inner QP through the \(T_k^{-r}\) and \(T_k^{-(2r-1)}\) scaling factors.

To solve the outer problem, we use L-BFGS \cite{byrd1995limited} to compute a search direction, followed by an Armijo backtracking line search \cite{nocedal2006numerical}. If the L-BFGS direction fails to yield a valid descent direction, the solver falls back to standard steepest descent.

\subsection{Experimental design}
The sweep uses the following problem and optimizer settings:
\begin{table}[h]
\centering
\label{tab:experimental_design}
\begin{tabular}{ll}
\toprule
\textbf{Quantity} & \textbf{Value} \\
\midrule
Bézier segment degree $n$ & 4 \\
Continuity order $c$ & 2 \\
Smoothness weights & $[0.1, 0.1, 1.0, 0.1]$ \\
Boundary conditions & Zero initial/terminal velocity \& acceleration \\
Derivative bounds & Specified for orders 1 and 2 \\
Initial duration $T_k$ & 3.0 \\
Time penalty $\lambda$ & 100 \\
L-BFGS iterations & 35 \\
L-BFGS memory & 10 \\
Armijo constant $c_1$ & $10^{-4}$ \\
Backtracking factor & 0.5 \\
Maximum backtracking steps & 30 \\
Duration floor $T_{\min}$ & 0.1 \\
Step stopping tolerance & $10^{-8}$ \\
\bottomrule
\end{tabular}
\end{table}

The tolerance-relaxation pairs in the sweep are
\begin{align}
&\texttt{solver\_tol}=10^{-4},\\
&\kappa_\mathrm{relax}=\{10^{0}, 10^{-1}, 10^{-2}, 10^{-3}, 10^{-4}\}.
\end{align}

\subsection{Complementary Results}
In addition to the results in Fig.~\ref{fig:forest}, we provide additional qualitative visualizations of the bilevel trajectory optimization operating in an office environment in Fig.~\ref{fig:office}.

\section{Learning a Safety Filter}\label{appendix:safety_filter}
\subsection{Problem formulation}
We consider $n$ planar agents with positions $p_i \in \mathbb{R}^2$ and stacked state and control
\begin{equation}
    x = \begin{bmatrix} p_1^\top & \cdots & p_n^\top \end{bmatrix}^\top \in \mathbb{R}^{2n},
    u = \begin{bmatrix} u_1^\top & \cdots & u_n^\top \end{bmatrix}^\top \in \mathbb{R}^{2n},
    \qquad
    p_{i,t+1} = p_{i,t} + \Delta t\, u_{i,t}.
\end{equation}
Each agent must reach a prescribed goal while avoiding static obstacles and pairwise collisions. The learned controller predicts a nominal control $u^{\mathrm{nom}}(x)$ together with state-dependent barrier gains $\alpha(x)$. These quantities parameterize a QP whose solution is the executed control $u^{\mathrm{safe}}(x)$. The neural network therefore captures task intent, while the QP acts as a safety-preserving projection that modifies the nominal action only when constraints become active.

\subsection{Safety filter as a CBF}
Let $m$ denote the number of circular obstacles, with center $c_o$ and radius $r_o$. Obstacle avoidance is encoded by
\begin{equation}
    h_{i,o}^{\mathrm{obs}}(x) = \|p_i - c_o\|_2^2 - r_o^2.
\end{equation}
For pairwise safety, with agent radius $r_{\mathrm{agent}}$,
\begin{equation}
    h_{i,j}^{\mathrm{pair}}(x) = \|p_i - p_j\|_2^2 - (2r_{\mathrm{agent}})^2.
\end{equation}
Under single-integrator dynamics, each barrier is enforced through the zeroing-CBF condition
\begin{equation}
    \nabla h(x)^\top u + \alpha(x)\, h(x) \geq 0,
\end{equation}
with learned positive gain $\alpha(x)$. The safety filter is the centralized QP
\begin{align}
    u^{\mathrm{safe}}(x)
    = \arg\min_{u \in \mathbb{R}^{2n}} \;& \frac{1}{2}\|u - u^{\mathrm{nom}}(x)\|_2^2 \label{eq:sf_qp_obj} \\
    \text{subject to} \;& \nabla_{p_i} h_{i,o}^{\mathrm{obs}}(x)^\top u_i + \alpha_{i,o}(x)\, h_{i,o}^{\mathrm{obs}}(x) \geq 0,
    \quad \forall i,o, \nonumber \\
    & \nabla h_{i,j}^{\mathrm{pair}}(x)^\top u + \alpha_{i,j}(x)\, h_{i,j}^{\mathrm{pair}}(x) \geq 0,
    \quad \forall i<j, \nonumber \\
    & -u_{\max} \leq u_k \leq u_{\max}, \quad k=1,\dots,2n. \nonumber
\end{align}
This formulation is minimally invasive: it returns the closest safe control to the nominal proposal. Since pairwise constraints couple multiple agents, the QP is solved jointly over all controls.

\subsection{Data generation}
Training data are generated by rolling out an oracle safety filter over randomly sampled obstacle configurations. At each time step, an unconstrained goal-seeking controller proposes
\begin{equation}
    u_{i}^{\mathrm{pd}}(x)
    =
    k_p \frac{g_i - p_i}{\max(\|g_i - p_i\|_2,\varepsilon)},
\end{equation}
which is clipped to the control bounds and passed through the oracle CBF-QP with fixed barrier gain. The executed safe action is stored together with the current state and obstacle layout. This yields a supervised dataset of dynamically relevant states, rather than arbitrary random samples, with concentration near regions where safety constraints become active.

\subsection{Training}
Training is performed by imitation learning. Given oracle-safe actions $u_\star^{(k)}$, we minimize
\begin{equation}
    \mathcal{L}(\theta)
    =
    \frac{1}{N}
    \sum_{k=1}^{N}
    \left\|
        \pi_\theta(x^{(k)}, \mathcal{O}^{(k)}) - u_{\star}^{(k)}
    \right\|_2^2,
\end{equation}
where $\pi_\theta$ is the learned controller and $\mathcal{O}^{(k)}$ denotes the obstacle layout for the $k$-th sample.  We use Adam~\cite{kingma2014adam} with learning rate $10^{-3}$ for $50$ epochs and batch size $256$. In the nominal experiment, $\Delta t=0.1$, $u_{\max}=2.0$, the goal tolerance is $0.1$, there are $3$ obstacles, and $r_{\mathrm{agent}}=0.35$. The study uses $n=7$ agents, $400$ training environments, $40$ test environments, and rollouts of length up to $200$. The resulting centralized QP has $14$ variables and
\begin{equation}
    4n + nm + \frac{n(n-1)}{2}
    = 28 + 21 + 21
    = 70
\end{equation}
constraints. We compare two differentiable QP backends, implicit and explicit, each in float64 and float32 precision, together with a neural baseline without a QP layer.

\subsection{Evaluation/Validation}
Validation is performed by closed-loop rollout on unseen obstacle configurations. We measure obstacle violations, pairwise violations, and whether all agents reach their goals within the horizon. In the nominal seven-agent setting, the successful CBF-QP models achieve zero safety violations on held-out rollouts, whereas the unconstrained neural baseline exhibits frequent collisions and much weaker task completion. Among the barrier-based variants, float64 and float32 implicit and float64 explicit train and validate successfully, while float32 explicit is numerically unstable in this regime.

\subsection{Maximum Batch size analysis}
To assess scalability, we vary the number of agents from $n=2$ to $n=9$ and record the largest batch size that can be trained without out-of-memory failure. This study is carried out for implicit in both float64 and float32 precision. As expected, the feasible batch size decreases rapidly with $n$ because the centralized QP grows in both dimension and constraint count, with the pairwise term scaling quadratically in the number of agents. Across the tested range, float32 supports consistently larger batches than float64, by roughly a factor of two, revealing a clear tradeoff between numerical robustness and memory efficiency.

\subsection{Numerical Details}
Throughout the experiment, the solver-specific QP settings are:
\begin{align}
  \texttt{f64}:&\quad \texttt{solver\_tol}=10^{-8},\;\kappa_\mathrm{relax}=10^{-4}, \\
  \texttt{f32}:&\quad \texttt{solver\_tol}=10^{-3},\;\kappa_\mathrm{relax}=10^{-3}.
\end{align}

The specific configuration used for data generation, training and validation is summarized below:
\begin{table}[h]
  \centering
  \label{tab:nominal-config}
  \begin{tabular}{ll}
    \toprule
    Quantity & Value \\
    \midrule
    Time step $\Delta t$ & $0.1$ \\
    Control bound $u_{\max}$ & $2.0$ \\
    Goal tolerance during data generation & $0.1$ \\
    Number of agents $n$ & $1$ to $9$ \\
    Number of obstacles $m$ & $3$ \\
    Workspace & $[-1,11]^2$ \\
    Agent radius $r_{\mathrm{agent}}$ & $0.35$ \\
    Oracle solver & \texttt{explicit f64} \\
    Training obstacle configurations & $400$ \\
    Test obstacle configurations & $40$ \\
    Maximum rollout length & $200$ steps \\
    Fixed oracle barrier gain \texttt{alpha\_fixed} & $1.0$ \\
    PD gain $k_p$ & $2.5$ \\
    Training epochs & $50$ \\
    Learning rate & $10^{-3}$ \\
    Batch size for every nominal model & $256$ \\
    \bottomrule
  \end{tabular}
\end{table}

\subsection{Complementary results}
We show random configurations obtained from learning the multi-agent safety filter for various numbers of agents and different QP solver configurations:
\begin{figure}[h]
  \centering
  \includegraphics[width=0.85\linewidth]{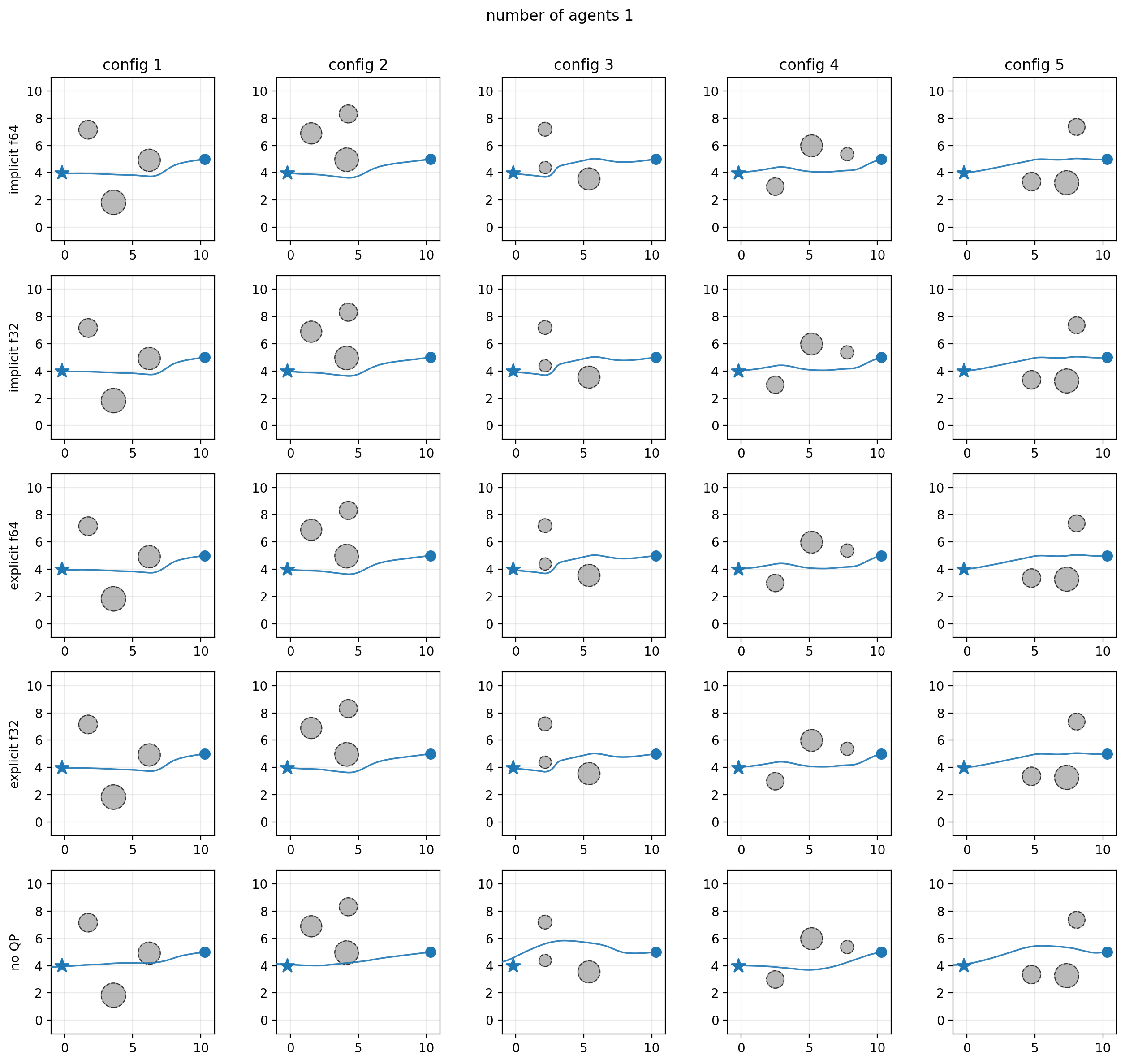}
  \caption{A subset of testing configurations for 1 agent. Only case where explicit f32 works.}
\end{figure}
\begin{figure}[h]
  \centering
  \includegraphics[width=0.9\linewidth]{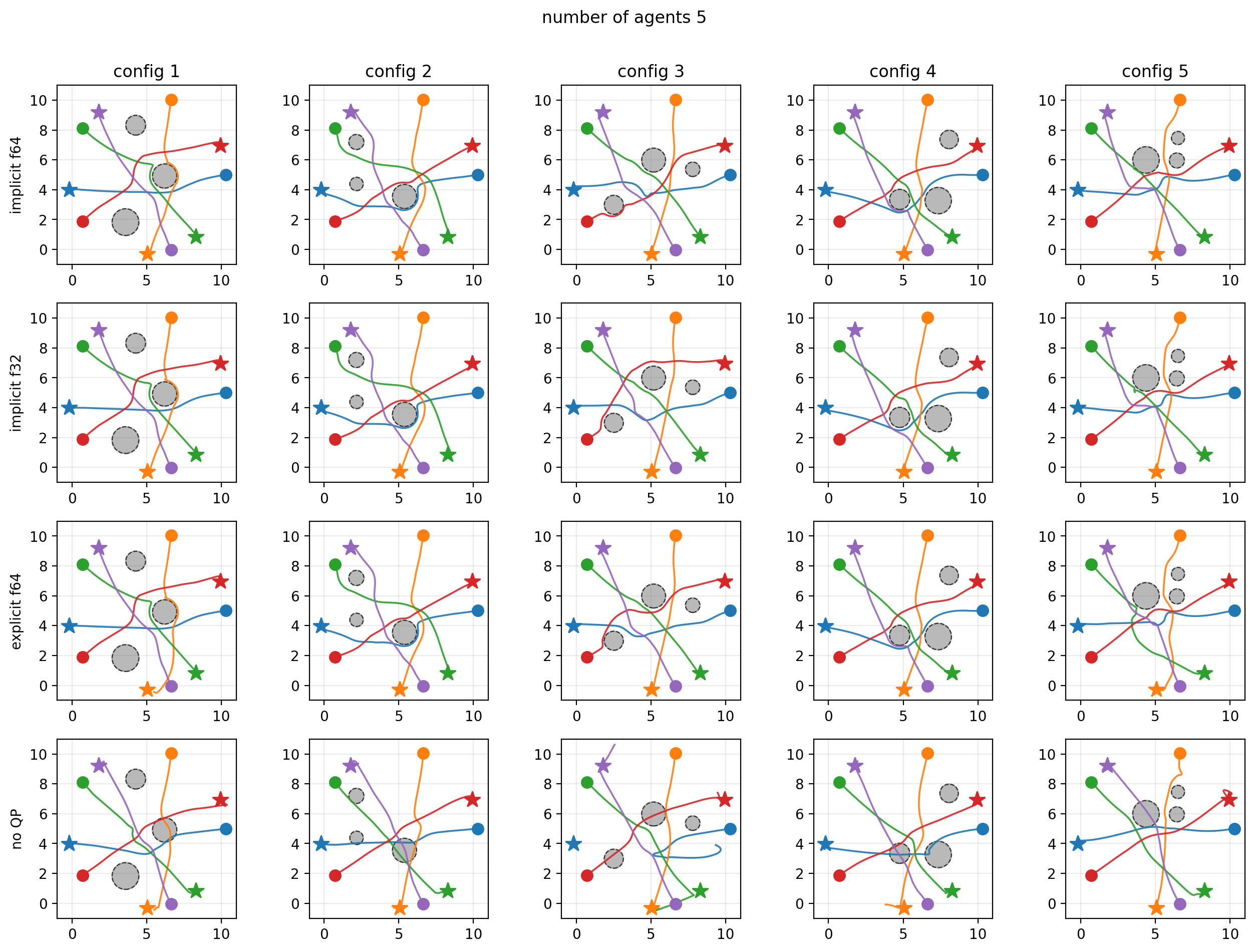}
  \caption{A subset of testing configurations for 3 agents.}
\end{figure}
\begin{figure}[h]
  \centering
  \includegraphics[width=0.95\linewidth]{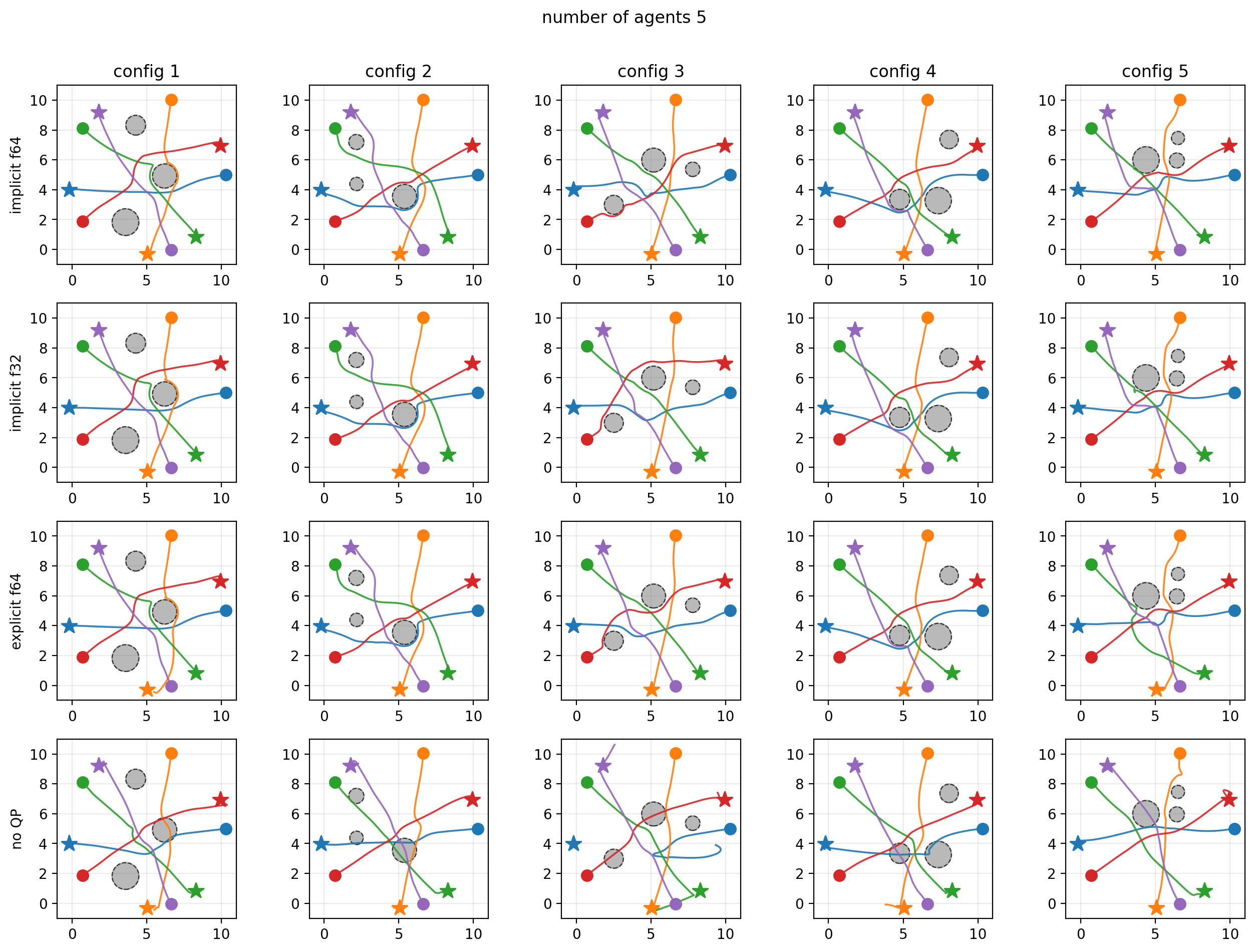}
  \caption{A subset of testing configurations for 5 agents.}
\end{figure}
\clearpage
\begin{figure}[t!]
  \centering
  \includegraphics[width=0.95\linewidth]{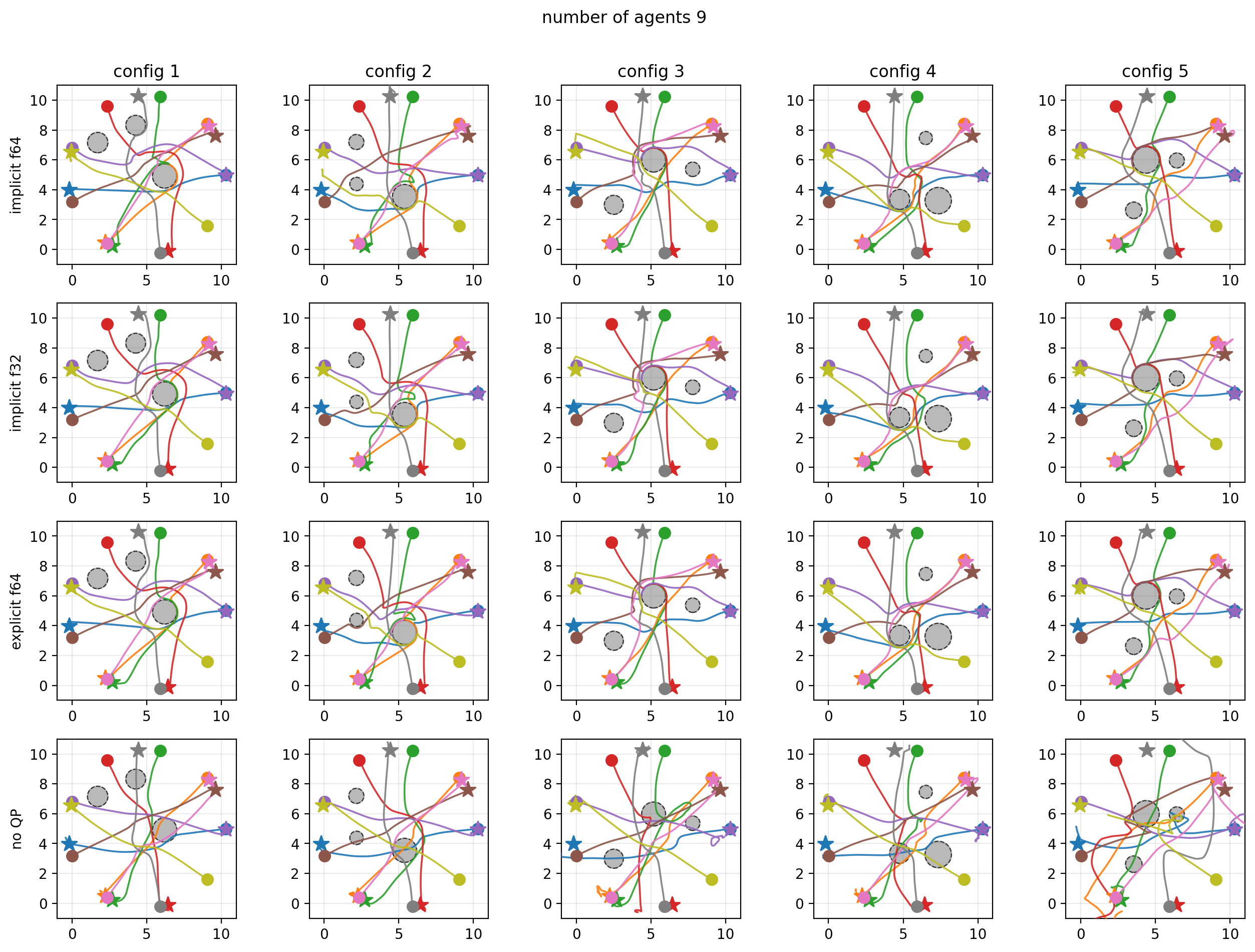}
  \caption{A subset of testing configurations for 9 agents.}
\end{figure}

\mbox{}

\ifarxiv
\fi

\end{document}